\documentclass[12pt]{article}

\usepackage[T2A]{fontenc}
\usepackage[utf8]{inputenc}

\usepackage{amsfonts}
\usepackage{amssymb}
\usepackage{amsmath}
\usepackage{amsthm}

\def \le {\leqslant}
\def \ge {\geqslant}

\theoremstyle{plain}

\begin{document}

\begin{Huge}
\centerline{\bf Exponents for three-dimensional} \centerline{\bf
simultaneous Diophantine approximations}
\end{Huge}
\begin{Large}
\vskip+0.5cm \centerline{by Nikolay  Moshchevitin\footnote{
Research is supported by RFBR grant No. 09-01-00371a  }}
\end{Large}
\vskip+2.0cm

\begin{small}
\centerline{\bf Abstract}
Let $\Theta = (\theta_1,\theta_2,\theta_3)\in \mathbb{R}^3$. Suppose that $1,\theta_1,\theta_2,\theta_3$ are linearly independent over $\mathbb{Z}$.
For Diophantine exponents
$$
\alpha(\Theta ) = \sup
\{\gamma >0:\,\,\,
\limsup_{t\to +\infty} t^\gamma \psi_\Theta (t) <+\infty \}
,$$
$$\beta(\Theta) = \sup
\{\gamma >0:\,\,\, \liminf_{t\to +\infty} t^\gamma \psi_\Theta (t)
<+\infty\}
$$
we prove
$$
\beta (\Theta) \ge
\frac{1}{2} \left( \frac{\alpha (\Theta)}{1-\alpha(\Theta)}
+\sqrt{\left(\frac{\alpha(\Theta)}{1-\alpha(\Theta)}\right)^2
+\frac{4\alpha(\Theta)}{1-\alpha(\Theta)}}\right) \alpha (\Theta)
$$


\end{small}

{\bf 1. Diophantine exponents.}

Let $\Theta =(\theta_1, ....,\theta_n)$ be a  real vector. We deal
with the function
$$
\psi_\Theta (t) = \min_{x\le t} \max_{1\le i\le n} ||\theta_i x||
.
$$
Here the minimum is taken over positive integers $x$ and $||\cdot
||$ stands for the distance to the nearest integer.

Suppose that at least one of the numbers $\theta_1, ....,\theta_n$
is irrational. Then
  $\psi_\Theta ( t) >0$ for all $t\ge 1$.  The
{\it uniform} Diophantine exponent $\alpha (\Theta)$ is defined as
the supremum of the set
$$\{\gamma >0:\,\,\,
\limsup_{t\to +\infty} t^\gamma \psi_\Theta (t) <+\infty \},
$$
It is a well-known fact that for all $\Theta$ one has
$$
\frac{1}{n}\le\alpha (\Theta)\le 1.
$$
 The
{\it ordinary} Diophantine exponent $\beta (\Theta)$ is defined as
the supremum of the set
$$
\{\gamma >0:\,\,\, \liminf_{t\to +\infty} t^\gamma \psi_\Theta (t)
<+\infty\} .
$$
Obviously
\begin{equation}\label{ab}
\beta (\Theta) \ge \alpha (\Theta).
\end{equation}

 {\bf 2. Functions.}

 For each $ \alpha \in \left[\frac{1}{3}, 1\right)$,
define
$$
g_1 (\alpha ) = \frac{\alpha}{1-\alpha}
$$
and
$$
g_2(\alpha ) = \frac{\alpha (1-\alpha) +\sqrt{
\alpha(\alpha^3+6\alpha^2-7\alpha+4)}}{2(2\alpha^2-2\alpha+1)}.
$$
The value $g_2(\alpha)$ is the largest root of the equation
$$
(2\alpha^2-2\alpha +1)x^2+\alpha(\alpha -1)x -\alpha = 0.$$ Note
that
$$
g_2(1/3)=g_2(1)=1,
$$
and for $1/3<\alpha <1$ one has $g_2(\alpha )>1$. Let $\alpha_0$
be the unique real root of the equation
$$
x^3-x^2+2x - 1 =0.
$$
In the interval $ 1/3< \alpha <\alpha_0$ one has
\begin{equation}\label{staro} g_2(\alpha ) > \max \left(1, g_1(\alpha
)\right).
\end{equation}
In the interval $\alpha_0 \le \alpha < 1$ we see that
$$
g_2(\alpha ) \le g_1(\alpha).
$$
We define one more function. Put
\begin{equation}\label{functions}
g_3(\alpha ) = \frac{1}{2} \left( \frac{\alpha}{1-\alpha}
+\sqrt{\left(\frac{\alpha}{1-\alpha}\right)^2
+\frac{4\alpha}{1-\alpha}}\right) .
\end{equation}
 Simple calculation shows that

\begin{equation}\label{bet1}
g_3(\alpha)  > \max (g_1 (\alpha ), g_2(\alpha ) )\,\,\,\,\,
\forall \, \alpha \in \left( \frac{1}{3}, 1 \right).
\end{equation}

{\bf 3. Jarn\'{\i}k's result.}

In a fundamental paper \cite{J} V. Jarn\'{\i}k proved the following
theorem.

{ \bf Theorem 1.}

{\it  Let $\psi(t)$ be a continuous function in $t$, decreasing to
zero as $t\to+\infty$. Suppose that the function $t\psi(t)$
increases to infinity as
 $t\to +\infty$.
Let~$\rho(t)$ be the inverse function to the function ~$t\psi(t)$.
Put
$$
\varphi^{[\psi]}(t)=\psi\biggl(\rho\biggl(\frac{1}{6\psi(t)}\biggr)\biggr).
$$
Suppose that $n \ge 2$ and  among numbers $\theta_1,...,\theta_n$
there exist at least two numbers which, together with 1, are
linearly independent over $\mathbb{Z}$. Suppose that
$$
\psi_\Theta (t) \le \psi (t) $$ for all $t$ large enough. Then
there exist infinitely many integers $x$ such that
$$
\max_{1\le j \le n}\|x\theta_j\|\le \varphi^{[\psi ]}(x).
$$
}
 The next Jarn\'{\i}k's result on Diophantine exponents is an
obvious corollary of Theorem 1.
\newpage
{ \bf Theorem 2.}

{\it Suppose that $n \ge 2$ and  among numbers
$\theta_1,...,\theta_n$ there exist at least two numbers which,
together with 1, are linearly independent over $\mathbb{Z}$. Then
$$
\beta (\Theta ) \ge \alpha (\theta) g_1(\alpha (\Theta)).
$$}

To obtain Theorem 2 from Theorem 1 one takes $\psi (t) =
t^{-\alpha} $ with $ \alpha < \alpha (\Theta)$.

From the other hand V. Jarn\'{\i}k  \cite{J}  proved that there
exists a collection of numbers $\Theta = (\theta_1,...,\theta_n)$
such that $1,\theta_1,...,\theta_n$ are linearly independent over
$\mathbb{Z}$  and
$$
\beta(\Theta )< \frac{\alpha(\Theta)}{1-\alpha(\Theta)}.
$$

 In the case $n=2$ the lower bound of
Jarn\'{\i}k's Theorem 2 is optimal. The following result was
proved by M. Laurent \cite{Lo}.

{\bf Theorem 3.}

{\it
For any $\alpha, \beta >0$  satisfying
$$
\frac{1}{2}\le \alpha \le 1,\,\,\,\, \beta \ge \alpha g_1(\alpha )
$$
there exists a vector $\Theta = (\theta_1,\theta_2) \in \mathbb{R}^2$ such that
 $$
 \alpha (\Theta ) =\alpha,\,\,\,\, \beta (\Theta ) = \beta.
 $$
}

This result is a corollary of a general theorem concerning four two-dimensional Diophantine exponents.

Note that in the case $n \ge 3$ the bound of
Theorem 2 in the range $\frac{1}{n}\le \alpha < \frac{1}{2}$ is
weaker than the trivial bound (\ref{ab}).

 N. Moshchevitin \cite{Mo1} (see also \cite{Mo2}, Section 5.2)  improved Jarn\'{\i}k's
result in the case $n=3$ and for $\alpha \in
\left(\frac{1}{3},{\alpha_0}\right)$. He obtained the following

 {\bf Theorem 4.}\,\,{\it Suppose that $m=1,n=3$
and the collection $ \Theta=\left(  \theta_1, \theta_2 , \theta_3
\right)$ consists of numbers which, together with 1, are linearly
independent over $\mathbb{Z}$. Then
$$
\beta (\Theta) \ge \alpha (\Theta) g_2 (\alpha (\Theta)).
$$
}

In the case $n=3$, Theorems 2 and 4 together give an estimate
which is better than the trivial estimate (\ref{ab}) for all
admissible values of $\alpha (\Theta)$.

{\bf 4. New result.}

In this paper we give a new lower bound for $\beta (\Theta)$ in
terms of $\alpha (\Theta)$. From (\ref{bet1}) it follows that this
bound is stronger than all previous bounds (Theorems 2 and 4) for
all admissible values of $\alpha (\Theta)$.

 {\bf Theorem 5.}\,\,

 {\it Suppose that $m=1,n=3$
and the vector $ \Theta=\left( \theta_1,\theta_2, \theta_3\right)$
consists of numbers linearly independent, together with 1, over
$\mathbb{Z}$. Then
$$
\beta (\Theta) \ge \alpha (\Theta) g_3 (\alpha (\Theta)).
$$
 }

Sections 5,6,7 below contains auxiliary results. Theorem 5 is
proved in Section 8.

\newpage
{\bf 5. Best approximations.}

For  each integer  $ { x} $, put
$$
  \zeta ({ x}) =  \max_{1\le j\le n}||\theta_j{ x}||.
$$
A positive integer
  $  x $ is said to be  a {\it best approximation} if
$$
\zeta ({ x})=\min_{{ x}'} \zeta ({ x}'),$$ where the minimum is
taken over all  $ { x}'  \in \mathbb{Z} $ such that
$$
0<  { x}'\le { x}. $$
 Consider the case when all numbers
1 and $\theta_{j},\,\,\,  1\le j\le n$  are linearly independent
over  $\mathbb{Z}$. Then all best approximations lead to sequences
  $${ x}_1< { x}_2<...<{ x}_\nu<{ x}_{\nu+1}<...
\,\,\, , $$
$$
\zeta({x}_1)> \zeta({ x}_2)>...>\zeta({ x}_\nu)>\zeta({
x}_{\nu+1})>... \,\,\, . $$ We use the notation
$$
 \zeta_\nu = \zeta({ x}_\nu).
$$
Choose $ y_{1,\nu},...,y_{n,\nu}\in \mathbb{Z}$ so that
$$
||\theta_j{ x}_\nu||= |\theta_j{\bf x}_\nu -y_{j,\nu}|.
$$
 We define
$$
{{\bf z}_\nu} = (x_{\nu}, y_{1,\nu},..., y_{n,\nu} ) \in
\mathbb{Z}^{n+1}.
$$
If $\psi (t)$ is a continuous function decreasing to $0$  as $
t\to \infty$, with
$$
\psi_{\Theta}(t) \le \psi (t)$$ then one  easily sees that
\begin{equation}\label{geret}
\zeta_\nu \le \psi (x_{\nu+1}).\end{equation}

 Some useful fact about best approximations can be  found in \cite{Mo2}.

{\bf 6. Two-dimensional subspaces.}

{\bf Lemma 1.}\,\,{\it Suppose that all vectors of the best
approximations $ {\bf z}_l , \nu\le l \le k$ lie in a certain
two-dimensional linear subspace $\pi \subset \mathbb{R}^4$.
Consider two-dimensional lattice $\Lambda = \pi \cap \mathbb{Z}^4$
with two-dimensional fundamental volume  $ {\rm det}\, \Lambda$.
 Then
for all $l$ from the interval $\nu\le l\le k-1$ one has
\begin{equation}\label{le1}
C_1 \, {{\rm det}\,\Lambda}  \le \zeta_{l}x_{l+1} \le 2\,{\rm det}
\Lambda.
\end{equation}
where $C_1 =
\left(2\sqrt{3\left(1+\left(|\theta_1|+\frac{1}{2}\right)^2+\left(|\theta_2|+\frac{1}{2}\right)^2+
\left(|\theta_3|+\frac{1}{2}\right)^2\right)}\right)^{-1}$.
 In particular,
\begin{equation}\label{le2}
{\rm det}\, \Lambda\ge \frac{ \min (\zeta_{\nu}x_{\nu+1} ,
\zeta_{k-1}x_k)}{2}.
\end{equation}
}

{\bf  Proof.}\,\,  The parallelepiped
$$
\Omega_l = \{ {\bf z}=(x,y_1,y_2,y_3):\,\, |x|< x_{l+1},\,\,
\max_{1\le j\le 3}|\theta_jx-y_j|<\zeta_l\}
$$
has no non-zero integer points inside for every $l$. Consider
two-dimensional ${\bf 0}$-symmetric convex body
$$
\Xi_l =\Omega_l \cap \pi.
$$
  One can see that the two-dimensional  Lebesgue measure $\mu (\Xi_l)$ of
$\Xi_l$ admit the following lower and upper bounds:
\begin{equation}\label{mea}
2\zeta_l x_{l+1} \le \mu (\Xi_l)\le
4\sqrt{3\left(1+\left(|\theta_1|+\frac{1}{2}\right)^2+\left(|\theta_2|+\frac{1}{2}\right)^2+
\left(|\theta_3|+\frac{1}{2}\right)^2\right)} \,\,\zeta_l x_{l+1}.
\end{equation}
We see that there is no non-zero points of $\Lambda$ inside
$\Xi_l$ and that there are two linearly independent points ${\bf
z}_l, {\bf z}_{l+1}\in \Lambda$ on the boundary of $\Xi_l$. So
obviously
\begin{equation}\label{mea1}
2\, {\rm det}\Lambda \le \mu (\Xi_l).
\end{equation}
From the Minkowski convex body theorem it follows that
\begin{equation}\label{mea2}
  \mu (\Xi_l)\le 4  \,{\rm det}\Lambda .
\end{equation}
 Now (\ref{le1}) follows from ({\ref{mea},\ref{mea1},\ref{mea2}). Lemma  is proved.$\Box$

{\bf 7. Three-dimensional subspaces.}

Consider three consecutive best approximation vectors ${\bf
z}_{l-1}, {\bf z}_l, {\bf z}_{l+1}$. Suppose that these vectors
are linearly independent. Consider the three-dimensional linear
subspace $$ \Pi_l ={\rm span} ({\bf z}_{l-1}, {\bf z}_l, {\bf
z}_{l+1}).$$ Consider the lattice
$$
\Gamma_l = \Pi_l \cap \mathbb{Z}^4
$$
with the fundamental volume ${\rm det} \, \Gamma_l$. Let $\Delta$
be three-dimensional volume of the three-dimensional simplex
${\cal S}$ with vertices ${\bf 0},{\bf z}_{l-1}, {\bf z}_l, {\bf
z}_{l+1}$. We see that
\begin{equation}\label{nona}
\Delta \ge \frac{{\rm det} \, \Gamma_l}{6}.
\end{equation}
Consider determinants
\begin{equation}\label{determa}
\Delta_1= -\left| \begin{array}{ccc} x_{l-1} & y_{2,l-1} &
y_{3,l-1} \cr x_{l} & y_{2,l} & y_{3,l} \cr x_{l+1} & y_{2,l+1} &
y_{3,l+1}
\end{array}
\right|, \,\,\Delta_2=  \left| \begin{array}{ccc} x_{l-1} &
y_{1,l-1} & y_{3,l-1} \cr x_{l} & y_{1,l} & y_{3,l} \cr x_{l+1} &
y_{1,l+1} & y_{3,l+1}
\end{array}
\right|, \,\,
 \Delta_3= -\left| \begin{array}{ccc} x_{l-1}
& y_{1,l-1} & y_{2,l-1} \cr x_{l} & y_{1,l} & y_{2,l} \cr x_{l+1}
& y_{1,l+1} & y_{2,l+1}
\end{array}
\right|.
\end{equation}
Absolute values of these determinants are equal to
three-dimensional volumes of projections of the simplex ${\cal S}$
onto three-dimensional coordinate subspaces ($\{y_1 =0\}, \{y_2 =
0\}$ and $\{ y_3=0\}$ respectively) multiplied by $6$.

Note that for $j =1,2,3$ one has
\begin{equation}\label{star}
|\Delta_j| \le 6\zeta_{l-1}\zeta_l x_{l+1}.
\end{equation}

 {\bf Lemma 2.}\,\,{\it Among determinants (\ref{determa}) there
exist a determinant with absolute value $\ge C_2 \Delta$, where
$C_2 = 2/
 (2+ \max_{1\le i\le 3} |\theta_i|)$.}

{\bf Proof.}

Consider the determinant
$$
 \Delta_0= \left| \begin{array}{ccc} y_{1,l-1} &
y_{2,l-1} & y_{3,l-1} \cr y_{1,l} & y_{2,l} & y_{3,l} \cr
y_{1,l+1} & y_{2,l+1} & y_{3,l+1}
\end{array}
\right|
$$
and the vector
$$
{\bf w} = (\Delta_0, \Delta_1, \Delta_2, \Delta_3) \in
\mathbb{Z}^4.
$$
We see that ${\bf w}$ is orthogonal to the subspace $\Pi_l$, that
is
$$
\Delta_0x_j+ \Delta_1y_{1,j}+ \Delta_2y_{2,j}+ \Delta_3y_{3,j} =
0,\,\,\,\,\, j = l-1,l,l+1.
$$
So
$$
\Delta_0 = -\sum_{i=1}^3 \Delta_i \frac{y_{i,l}}{x_l}=
-\sum_{i=1}^3 \Delta_i \left(\frac{y_{i,l}}{x_l}-\theta_i\right)
-\sum_{i=1}^3 \Delta_i \theta_i.
$$
As $\left|\frac{y_{i,l}}{x_l}-\theta_i\right| \le 1$ we see that
\begin{equation}\label{trivi}
|\Delta_0| \le ( 1+\max_{1\le i\le 3} |\theta_i|)
(|\Delta_1|+|\Delta_2|+|\Delta_3|).
\end{equation}
But
\begin{equation}\label{summ}
36 \Delta^2 = \Delta_0^2+ \Delta_1^2+\Delta_2^2+\Delta_3^2.
\end{equation}
 From (\ref{trivi},\ref{summ}) we deduce the inequality
 $$
\Delta \le \frac{1}{6} \,(2+ \max_{1\le i\le 3} |\theta_i|)\,
(|\Delta_1|+|\Delta_2|+|\Delta_3|),
 $$
and the lemma follows.$\Box$

{\bf 8. Proof of Theorem 5.}

Take $\alpha < \alpha(\Theta)$. Then
\begin{equation}\label{www}
\zeta_l \le x_{l+1}^{-\alpha}
\end{equation} for all $l$
large
enough.

 Consider best approximation vectors ${\bf z}_\nu =
(x_\nu,y_{1,\nu},y_{2,\nu},y_{3,\nu})$. From the condition that
numbers $1,\theta_1,\theta_2,\theta_3$ are linearly independent
over $\mathbb{Z}$ we see that
   there exist infinitely many pairs of indices
   $\nu<k, \nu\to +\infty$   such that

$\bullet $  both triples
$$
{\bf z}_{\nu-1},{\bf z}_\nu,{\bf z}_{\nu+1};\,\,\,\,\,\, {\bf
z}_{k-1},{\bf z}_k,{\bf z}_{k+1}
$$ consist of
  linearly independent vectors;

$\bullet$ there exists a two-dimensional linear subspace   $\pi$
such that
$$
{\bf z}_l\in \pi,\,\,\, \nu\le l\le k;\,\,\,\,\, {\bf z}_{\nu-1}
\not\in \pi,\,\,\, {\bf z}_{k+1} \not\in \pi;
$$

$\bullet$  the vectors
$$
{\bf z}_{\nu-1},{\bf z}_\nu,{\bf z}_{k},{\bf z}_{k+1}
$$
are linearly independent.

Consider the two-dimensional lattice $$ \Lambda =\pi \cap \mathbb{Z}^4
$$
By Lemma 1, its two-dimensional  fundamental volume ${\rm
det}\, \Lambda$ satisfies
\begin{equation}\label{lll}
{\rm det}\, \Lambda \asymp_\Theta \zeta_\nu x_{\nu+1}
\asymp_\Theta \zeta_{k-1}x_k.
\end{equation}
Consider the two dimensional orthogonal complement $\pi^\perp$ to
$\pi$ and the lattice $$\Lambda^\perp = \pi^\perp \cap
\mathbb{Z}^4.$$ It is well-known that
\begin{equation}\label{lll1}
{\rm det}\, \Lambda^\perp = {\rm det}\, \Lambda.
\end{equation}
Consider the lattices
$$
\Gamma_\nu = ({\rm span}\, ( {\bf z}_{\nu-1},{\bf z}_\nu,{\bf
z}_{\nu+1})) \cap \mathbb{Z}^4,\,\,\,\, \Gamma_k = ({\rm span}\, (
{\bf z}_{k-1},{\bf z}_k,{\bf z}_{k+1}))\cap \mathbb{Z}^4.
$$
and primitive integer vectors ${\bf w}_\nu, {\bf w}_k\in
\mathbb{Z}^4$ which are orthogonal to $\Pi_\nu = {\rm span}\, (
{\bf z}_{\nu-1},{\bf z}_\nu,{\bf z}_{\nu+1})$, $\Pi_k = {\rm
span}\, ( {\bf z}_{k-1},{\bf z}_k,{\bf z}_{k+1})$ respectively.
Obviously
$$
{\bf w}_\nu, {\bf w}_k\in \Lambda^\perp.
$$
Put
$$
b= \frac{1}{2} \left( -\frac{\alpha}{1-\alpha} + \sqrt{ \left(
\frac{\alpha}{1-\alpha}\right)^2+ \frac{4\alpha}{1-\alpha}}\right)
\in (0,1), \,\,\,\,\, a = 1 - b,
$$
so
$$
\frac{\alpha}{1-\alpha}+ b =g_3(\alpha ).
$$
Then
$$
{\rm det}\, \Lambda^\perp \le |w_\nu|\cdot |w_k|,
$$
where $|\cdot |$ stands for the Euclidean norm, and so we obtain
that either
\begin{equation}\label{AA}
{\rm det }\, \Gamma_\nu = |{\bf w}_\nu| \ge ({\rm det
}\,\Lambda^\perp)^a = ({\rm det }\,\Lambda )^a
\end{equation}
or
\begin{equation}\label{BB}
{\rm det }\, \Gamma_k = |{\bf w}_k| \ge ({\rm det
}\,\Lambda^\perp)^b = ({\rm det }\,\Lambda )^b
\end{equation}
(using (\ref{lll1})).

If (\ref{AA}) holds then by Lemma 2, (\ref{star}), (\ref{nona})
and (\ref{lll})
  we see that
$$
\zeta_{\nu -1} \zeta_\nu x_{\nu+1} \gg |\Delta_j| \gg_\Theta {\rm
det }\, \Gamma_\nu  \gg_\Theta  ({\rm det }\,\Lambda )^a\gg
(\zeta_\nu x_{\nu+1})^a
$$
(here $\Delta_j$ is the determinant from Lemma 2 applied to the
lattice $\Gamma = \Gamma_\nu$). From the definition of $a$ and
(\ref{www}) we see that
$$
x_{\nu+1} \gg_\Theta x_\nu ^{g_3(\alpha)}.
$$
We apply (\ref{www}) again to obtain
$$
\zeta_\nu \ll_\Theta  x_\nu^{-\alpha g_3(\alpha)}.
$$

If (\ref{BB}) holds then by Lemma 2, (\ref{star}), (\ref{nona}) and (\ref{lll})
  we see that
$$
\zeta_{k -1} \zeta_k x_{k+1} \gg |\Delta_{j'}| \gg_\Theta {\rm det
}\, \Gamma_k  \gg_\Theta  ({\rm det }\,\Lambda )^b\gg (\zeta_{k-1}
x_{k})^b
$$
(here $\Delta_{j'}$ is the determinant from Lemma 2 applied to the
lattice $\Gamma = \Gamma_k$). From the definition of $b$ and
(\ref{www}) we see that
$$
x_{k+1} \gg_\Theta x_k ^{g_3(\alpha)}.
$$
We apply (\ref{www}) again to obtain
$$
\zeta_k \ll_\Theta  x_k^{-\alpha g_3(\alpha)}.
$$
Theorem 5 is proved.$\Box$
\newpage
{\bf 8. Acknowledgement and a remark.}

The author thanks the anonymous referee for useful and important
suggestions. Here we would like to note that the referee pointed
out that it is possible to get a simpler proof of Theorem 5 by
means of W.M. Schmidt's inequality on heights of rational
subspaces (see \cite{SCH67}). For a rational subspace $U\subset
\mathbb{R}^n$ its height $H(U)$ is defined as the co-volume of the
lattice $U\cap\mathbb{Z}^n$. Schmidt shows that for any two
rational subspaces $U,V \in \mathbb{R}^n$ one has
$$
H(U\cap V) H(U+V) \ll_n H(U) H(V),
$$
To prove our Theorem 5 one can use this inequality for
$$U={\rm span } ( {\bf z}_{\nu-1}, {\bf z}_\nu),\,\,\,\,
V = {\rm span } ( {\bf z}_{\nu}, {\bf z}_{\nu+1})
$$
and for
$$
U' = V = {\rm span } ( {\bf z}_{k-1}, {\bf z}_k),\,\,\, V' = {\rm span } ( {\bf z}_{k}, {\bf z}_{k+1}).
$$


\begin{thebibliography}{100}




\bibitem{J}
V. Jarn\'{\i}k,\,\, Contribution  \`{a} la th\'{e}orie des
approximations diophantiennes lin\'{e}aires et homog\`{e}nes,
Czechoslovak Math. J. 4 (1954), 330 - 353 (in Russian, French
summary).


\bibitem{Lo} M. Laurent,\,\,\,
Exponents of Diophantine approximations in dimension two.//
Canad.J.Math. 61, 1 (2009),165 - 189;  preprint available at
arXiv:math/0611352v1 (2006).

\bibitem{Mo1}
N.G. Moshchevitin,\,\,\ Contribution to Vojt\v{e}ch
Jarn\'{\i}k.//Preprint available at arXiv:0912.2442v3 (2009).

\bibitem{Mo2}
N.G. Moshchevitin,\,\,\, Khintchine's singular Diophantine systems and their
applications.// Russian Mathematical Surveys. 65:3 43 - 126 (2010); Preprint
available at arXiv:0912.4503v1 (2009).

\bibitem{SCH67}
W.M. Schmidt,\,\,\, On
 heights of algebraic subspaces and Diophantine approximations. //
Ann. of Math. (2) 85 (1967), 430 - 472.


\end{thebibliography}
\end{document}